\documentclass[a4paper,12pt]{amsart}
\pdfoutput=1 
\usepackage{amssymb} 
\usepackage{amsfonts}
\usepackage{amsmath, amscd}
\usepackage{amsthm}
\usepackage[utf8]{inputenc}
\usepackage{pgf}
\usepackage{graphicx}
\usepackage{tikz}
\usepackage{pdfpages}
\usepackage[all]{xy}

\newcommand{\ctor}{(\C^\times)^2}
\newcommand{\rtor}{(\R^\times)^2}

\usetikzlibrary{patterns}
\title{Tropical limit of log-inflection points for planar curves}

\author{Grigory Mikhalkin}
\address{Grigory Mikhalkin, Université de Genève, Mathématiques, Villa Battelle, 1227 Carouge, Suisse}
\email{grigory.mikhalkin@unige.ch}

\author{Arthur Renaudineau}
\address{Arthur Renaudineau, Eberhard Karls Universität Tübingen, Fachbereich Mathematik, Auf der Morgenstelle 10,
72076 Tübingen, Deutschland}
\email{arren@math.uni-tuebingen.de}
\date{}

\newtheorem{theorem}{Theorem}
\newtheorem{example}[theorem]{Example}

\newtheorem{remark}[theorem]{Remark}
\newtheorem{corollary}[theorem]{Corollary}
\newtheorem{definition}[theorem]{Definition}
\newtheorem{defn}[theorem]{Definition}
\newtheorem{definition-proposition}[theorem]{Definition-Proposition}

\newtheorem{proposition}[theorem]{Proposition}

\newcommand{\R}{\mathbb{R}}

\newcommand{\Z}{\mathbb{Z}}
\newcommand{\CP}{\mathbb{CP}}
\newcommand{\QP}{\mathbb{QP}}

\newcommand{\TP}{\mathbb{TP}}

\newcommand{\C}{\mathbb{C}}
\newcommand{\T}{\mathbb{T}}
\newcommand{\trop}{\operatorname{trop}}
\begin{document}
\begin{abstract}
The paper
describes the behavior of
log-inflection points 
(i.e. points of
inflection with respect
to the parallelization of $\ctor$ given
by the multiplicative group law)
of curves in $\ctor$
under passing to the tropical limit.
Assuming that the limiting tropical curve is smooth,
we show that
log-inflection points accumulate by pairs at the midpoints
of bounded edges.

\end{abstract}
\maketitle
\section{Introduction}
\newcommand{\cpDelta}{\C P_\Delta}
\newcommand{\cp}{{\mathbb C}P}
\newcommand{\rp}{{\mathbb R}P}
\newcommand{\dd}{\partial}
\newcommand{\Log}{\mathrm{Log}}
\newcommand{\conj}{\mathrm{conj}}
Denote by $\C^\times=\C\setminus\{0\}$ the complex torus and
let $f:\ctor\to\C$
be a Laurent polynomial
in two variables with the Newton polygon
$\Delta$.
The zero locus
$$V_f=\{(z,w)\in\ctor\ |\ f(z,w)=0\}\subset\ctor$$
is a non-compact complex curve.
For non-singular $V_f$
the {\em logarithmic Gauss map} (see \cite{Ka-Gauss})
\begin{equation}
\gamma_f:V_f\to \cp^1
\end{equation}
is defined by $\gamma_f(z,w)=
(z\frac{\dd f}{\dd z}:
w\frac{\dd f}{\dd w})$.
This map depends only on $V_f$ and not on $f$
itself. It is the map that takes a point of $V_f$
to the slope of its tangent plane with respect
to the parallelization of $\ctor$ given by
the {\em multiplicative translations}
$(z,w)\mapsto (\alpha z,\beta w)$,
$\alpha,\beta\in\C^\times$. 
\begin{defn}
The critical points of $\gamma_f$ are called the
{\em log-inflection}
points of $V_f$.
\end{defn}

Denote the set of critical points 
by 
$\rho V_f\subset V_f$.
If all coefficients of the Laurent
polynomial $f$ are real
we define 
$$\R V_f=
\{(z,w)\in\rtor\ |\ f(z,w)=0\}\subset\rtor .$$
We have $$\gamma_f|_{\R V_f}:\R V_f\to\rp^1,$$
and we denote by $\rho \R V_f\subset \R V_f$ the set of critical points of $\gamma_f|_{\R V_f}$ (called real log-inflection points).
Recall the definition of tropical limits
 (see e.g. \cite{IKMZ}, \cite{BIMS} as well
as \cite{Mik06}). 
Let $$\{f_t=\sum\limits_{j,k\in\Z}
a_{jk,t}z^jw^k\}_{t\in A}$$ be a family of polynomials
parameterized by an (infinite) set $A$.
Suppose that the degree of $f_t$ is universally bounded
in $z$, $z^{-1}$, $w$ and $w^{-1}$.
Let $\alpha:A\to\R$ be a map whose
image is unbounded from above. 
The pair $(f_t,\alpha)$ is called
a {\em scaled sequence} of polynomials $f_t$
($\alpha$ is the scaling).

The scaled sequence of polynomials is said to
{\em converge tropically}
if the limit $$a_{jk}^{\trop}=
\lim\limits_{t\in A}\log_{\alpha(t)}|a_{jk,t}|
\ \ \in\R\cup\{-\infty\}$$
exists
for all $j,k\in\Z$. 
Here the limit is taken with
respect to the directed set $A$.
where the order is given by $\alpha$.
Define the limiting {\em Newton polygon}
$$\Delta=\operatorname{ConvexHull}
\{(j,k)\in\Z^2\ |\ a^{\trop}_{jk}\neq -\infty\}
\ \ \subset\R^2.$$
It is a bounded convex lattice polygon.
The set $\R\cup\{-\infty\}$ is called the
{\em set of tropical numbers} and denoted 
by
$\T$.
Denote by $\Log:\ctor\to\R^2$ the map defined by
$$\Log(z,w)=(\log|z|,\log|w|).$$
If a scaled sequence of polynomials
$(f_t,\alpha)$ converges tropically, then the limit
$$
C=\lim\limits_{t\in A}\mathrm{Log}_{\alpha(t)}(V_t)
$$ 
exists in the Hausdorff sense,
i.e. as the limit of closed subsets of $\R^2$
with the topology given by the Hausdorff
distance. It is called the {\em tropical limit}
of $\{V_t\subset\ctor\}_{t\in A}$.


More generally, given a scaled sequence $(V_t,\alpha)$
of curves $V_t\subset (\C^*)^n$ we say that it converges tropically
to $C\subset\R^n$ if $C$ is the limit of $\mathrm{Log}_{\alpha(t)}(V_t)$
in the Hausdorff sense.

It can be proved (see \cite{Mik05}, and \cite{IKMZ} for a more general statement)
that the limit $C$ is a rectilinear
graph in $\R^n$.
The edges $E$ of $C$ are straight intervals
with rational slopes, and can be
prescribed integer {\em weights} $w(E)$
(coming from $V_t$)
so that for every vertex $v\in C$ it holds the {\em balancing condition}
$$\sum\limits_{E\ni v} w(E)u(E)=0,$$
where $u(E)\in\Z^n$ is the primitive vector parallel to $E$
in the direction away from $v$.
A rectilinear graph $C\subset\R^n$ with these properties
is called a {\em tropical curve}. If $n=2$ then $C$
defines a lattice subdivision of 
the Newton polygon $\Delta$.
If each polygon of the subdivision
is a triangle of area $\frac 12$
(the minimal possible area for a lattice
polygon) then $C$ is called {\em smooth}.
For details we refer to \cite{IMS} and \cite{BIMS}.

\begin{definition}
Let $C$ be a smooth plane tropical curve. The {\em tropical parabolic locus} of $C$ is
the set $\rho C\subset C$ formed
by the midpoints of all bounded edges of $C$.
\end{definition}
The main result of this paper is the following.
\begin{theorem}
\label{parabolicforcurves}
Let $C\subset\R^2$
be a smooth tropical curve,
and  
$\{V_t\subset\ctor\}_{t\in A}$
be a family of complex curves
parameterised by the scaled sequence
$\alpha:A\to\R$.
Suppose that $V_t$ tropically converges to $C$
and that the limiting Newton polygon $\Delta$
coincides with the Newton polygons of $V_t$ for large $\alpha(t)$. 
Then
$\Log_{\alpha(t)}(\rho V_t)$
converges to $\rho C$ in the Hausdorff
metric when $\alpha(t)\to+\infty$.
Furthermore, a small neighborhood of 
a point from $\rho C$ contains exactly
two points of $\Log_{\alpha(t)}(\rho V_{t})$
for large $\alpha(t)$.
\end{theorem}

In other words, $2\rho C$ is the tropical
limit of $\rho V_t$. 
\begin{remark}
Theorem \ref{parabolicforcurves} provides a $2-1$ correspondence $\Phi_t:\rho V_t\to \rho C$ for large $\alpha(t)$. If the family $\{V_t\}_{t\in A}$ is defined over $\R$, then $\rho V_t$ is invariant under the involution $\conj$ of complex conjugation in $\ctor$. 	In particular, the pair $\Phi_t^{-1}(p)$ is $\conj$-invariant for each $p\in\rho C$.
\end{remark}

\textit{Acknowledgements.} Research is supported in part by the grants 159240, 159581, the NCCR SwissMAP project of the Swiss National Science Foundation and by DFG-grant MA 4797/6-1. The first author is thankful to Ernesto Lupercio for useful conversations
about tropical Gauss maps.

\section{The case of real curves}
Let $\{V_t\subset (\C^\times)^2\}_{t\in A}$ be a family of real curves 
(i.e. invariant with respect to the involution
of complex conjugation) tropically convergent to a smooth tropical curve
$C\subset\R^2$ with respect to a scaling $\alpha:A\to \R$. Assume that the limiting Newton polygon
coincides with the Newton polygons of $V_t$ for large $\alpha(t)$. 
Let $E$ be a bounded edge of $C$. Denote by $p\in E$ the midpoint of $E$,
and let $t\in A$ such that $\alpha(t)$ is so large that the 2-1 correspondence
$\Phi_t:\rho V_t\to \rho C$ is defined.
\begin{definition}
The edge $E$ is called {\em $V_t$-twisted} 
if $\conj$ preserves $\Phi_t^{-1}(p)$ pointwise. Otherwise $\conj$ induces a non-trivial permutation of $\Phi_t^{-1}(p)$, and we call $E$ {\em $V_t$-untwisted}. 
\end{definition}

Clearly $E$ is $V_t$-twisted if $\Phi_t^{-1}(p)\subset \R V_t$, and
$V_t$-untwisted if $\Phi_t^{-1}(p)$ forms a complex-conjugated pair. 
\begin{definition}{$($cf. \cite{BIMS}$)$}
\\A subset $T$ of the set of bounded edges of a smooth tropical curve
$C\subset\R^2$ is called
{\em twist-admissible} if they satisfy the following condition: 
for any cycle $\gamma$ of $C$ (considered as a graph), if $E_1,\cdots, E_k$ denote edges in $\gamma\cap T$, and if $u(E_i)$ is a primitive integer in the direction of $E_i$, then
\begin{equation}
\sum_{i=1}^k u(E_i)=(0,0) \mod 2.
\label{def:twistcondition}
\end{equation}

\end{definition}
Recall that {\em a patchworking
polynomial} may be defined as
$$f_t(z_1,z_2)=\sum\limits_{j\in B} \beta_j(t)
(\alpha(t))^{a_j}z^j.$$
Here $j=(j_1,j_2)$, $z^j=(z_1^{j_1},z_2^{j_2})$,
$j\mapsto a_j$
is a strictly convex function
from a finite subset $B\subset\Z^2$ to $\R$, 
and $\left(\beta_j(t)\right)_{t\in A}\subset\R$ is a family
converging (with respect to the directed set $A$) to a non-zero real number. 
\newcommand{\sign}{\operatorname{sign}}
Denote $$\sigma_j=\sign(\lim\limits_{t\in A} \beta_j(t))
\in \{-1,+1\}.$$
The family $\left(V_t\right)_{t\in A}$ tropically converge a tropical curve $C$ defined by the tropical polynomial with the coefficients $a_j$ (see e.g. \cite{IMS}). According to Viro's patchworking theorem
(cf. \cite{GKZ}, \cite{Viro}), if the tropical curve $C$ is smooth then the rigid isotopy class
of the curve
$V_t=\{(x,y)\in (\R^\times)^2\ |\ f_t(x,y)=0\}$
for large value of $\alpha(t)$ is determined by $C$
together with the signs $\sigma_j$ (see e.g. \cite{BIMS}).
\begin{proposition}
\label{prop:realinflectionpts}
Let $V_t$ be a family of real curves given by a patchworking polynomial as above, and let $C$ be its tropical limit. Assume that $C$ is smooth.
For each bounded edge $E$ of $C$, denote by $v_1^E$ and $v_2^E$ the two vertices of the segment $\Delta_E$ dual to $E$ in the subdivision of $\Delta$ defined by $C$. The segment $\Delta_E$ is adjacent to two other triangles of the subdivision of $\Delta$. Denote by $v_3^E$ and $v_4^E$ the vertices of these two triangles different from $v_1^E$ and $v_2^E$. Then $E$ is $V_t$-twisted if and only if
\begin{itemize}
\item $\sigma_{v_1^E}\sigma_{v_2^E}\sigma_{v_3^E}\sigma_{v_4^E}>0$ if the coordinates modulo $2$ of $v_3^E$ and $v_4^E$ are distinct,
\item $\sigma_{v_3^E}\sigma_{v_4^E}<0$ if the coordinates modulo $2$ of $v_3^E$ and $v_4^E$ are the same.
\end{itemize} 
\end{proposition}
%
The following corollary follows easily from Proposition \ref{prop:realinflectionpts} and known facts about patchworking (see \cite{BIMS}). 
\begin{corollary}
\label{corollary:twistededges}
The set of $V_t$-twisted edges of $C$ is twist-admissible.
Conversely, for any twist-admissible subset $T$ of bounded edges of
a smooth tropical curve $C$
there exist a family  $\{V_t\subset (\C^\times)^2\}_{t\in A}$
of real curves 
such that $T$ is the set of $V_t$-twisted edges.
\end{corollary}

\newcommand{\ignore}[1]{\relax}
\ignore{
We say that a family of real curve $\R V_t\subset\rtor$ (parameterized by a scaled sequence $\alpha:A\to \R$) tropically converges to $C\subset\R^2$ if the complexications $V_t\subset\ctor$ of $\R V_t$ tropically converge to $C$. Note that the group $(\Z/2\Z)^2$ acts on $\rtor$ via multiplication by $(\pm 1,\pm 1)$. Given $(\varepsilon_1,\varepsilon_2)\in (\Z/2\Z)^2$, we call the corresponding transformation $\sigma_{\varepsilon_1,\varepsilon_2}:\rtor\to\rtor$ a {\em coordinate sign changes}. 
\begin{theorem}
\label{realtheorem}
Let $C\subset\R^2$ be a smooth tropical curve, and $\{\R V_t\subset\rtor\}_{t\in A}$ be a family of real curves tropically converging to $C$. Then for large $\alpha(t_1),\alpha(t_2)>0$, the curve $\R V_{t_1}$ is isotopic to $\sigma_{\varepsilon_1,\varepsilon_2}(\R V_{t_2})$ in $\rtor$ for a certain coordinate sign change transformation $\sigma_{\varepsilon_1,\varepsilon_2}$ if and only if $V_{t_1}$ and $V_{t_2}$ have the same set of twisted edges on $C$. Furthermore, for any twist-admissible set $T$ of bounded edges of $C$ there exists a family of real curves tropically converging to $C$ and such that the set of twisted edges of real curves in this family close to the limit is $T$.
\end{theorem}
\begin{proof}
The first part follows immediatly from Theorem \ref{parabolicforcurves} and from Viro's patchworking theorem in its tropical reformulation (see for example \cite{BIMS} Section $3$). For the second part, choose a patchworking polynomial $f_t$ such that a bounded edge $E\notin T$ if and only if the amoeba of $V_t=\{f_t=0\}$ is smooth in a small neighborhood of $E$ for $t$ large enough (see \cite{BIMS} and \cite{Mik00} for a description of the stable singularities of the amoeba). Then the set of twisted edges of $C$ corresponding to $V_t$, for $t$ large enough, is exactly $T$.
\end{proof}
\begin{remark}
It follows from Theorem \ref{realtheorem} that phase-tropical convergence can be reinterpreted in terms of twisted edges. If $\{\R V_t\subset\rtor\}_{t\in A}$ is a family of real curves tropically converging to $C$, then it converges phase-tropically (up to the $(\Z/2\Z)^2$ action of coordinate sign changes) if $V_t$ has a constant set of twisted edges on $C$ for $\alpha(t)$ large enough. Moreover, if $\{\R V_t\subset\rtor\}_{t\in A}$ and $\{\R W_t\subset\rtor\}_{t\in B}$ are two families converging phase-tropically to $C$, then for $\alpha(t),\beta(t)$ large enough the curve $\R V_{t}$ is isotopic to $\sigma_{\varepsilon_1,\varepsilon_2}(\R W_{t})$ in $\rtor$ for a certain coordinate sign change transformation $\sigma_{\varepsilon_1,\varepsilon_2}$ if and only if $V_{t}$ and $W_{t}$ have the same set of twisted edges on $C$.
\end{remark}

\textit{Acknowledgements.} Research is supported in part by the grants 159240, 159581, the NCCR SwissMAP project of the Swiss National Science Foundation and by DFG-grant MA 4797/6-1. The first author is thankful to Ernesto Lupercio for useful conversations
about tropical Gauss maps.
}

\section{Tropical limit of logarithmic Gauss map}
\label{tropgaussmap}
Let $\{V_t\}_{t\in A}=\{f_t=0\}_{t\in A}\subset\ctor$
be a family of complex curves
parameterised by the scaled sequence
$\alpha:A\to\R$. Let
$\widetilde{V}_t$, $t\in A$, be the graph of the logarithmic Gauss map $\gamma_{f_t}$
(which we denote 
by
 $\gamma_t$ for brevity)
given by
$$
\widetilde{V}_t=\left\lbrace \left(z_t,\gamma_t(z_t)\right) \vert z_t\in V_t\right\rbrace \subset \ctor\times\CP^{1}.
$$ 
Denote by $\pi_1:\widetilde{V}_t\rightarrow V_t\subset\ctor$ the projection on the first factor and by $\pi_2:\widetilde{V}_t\rightarrow \CP^{1}$ the projection on the second factor. 
The map $\pi_1$ is an isomorphism and the map $\pi_2$ is of degree $2\mathrm{Area}(\Delta)$ by the Bernstein-Kouchnirenko formula (see \cite{Kou} and also \cite{Mik00}, Lemma $2$). 
Denote by $\TP^n$ the tropical projective space of dimension $n$ (see for example \cite{BIMS}). 
Denote by $\pi_1^{trop}:\R^2\times\TP^1\rightarrow \R^2$ the projection on the first factor
and by $\pi_2^{trop}:\R^2\times\TP^1\rightarrow \TP^{1}$ the projection on the second factor.  
\begin{proposition}
\label{compactnesstheorem}
There exists a subfamily $A'\subset A$ with
unbounded $\alpha(A')\subset\R$ such that
the subsequence $(\widetilde{V}_t)_{t\in A'}$ tropically converges to a tropical curve
$\widetilde{C}$ in $\R^2\times\TP^{1}$.
Furthermore, we have $\pi_1^{trop}(\widetilde{C})=C$.
\end{proposition}
\begin{proof}
The proposition is a special case of Theorem 38 of \cite{IKMZ}.
We apply this theorem to the projective curves $\bar V_t$, $t\in A$,
obtained as the closures of $\widetilde{V}_t\cap (\C^\times)^3$ in $\CP^3$.
We set $\widetilde{C}$ to be the closure in $\R^2\times\TP^1$ of
$\bar C\cap\R^3$, where $\bar C\subset\TP^3$ is the tropical limit
of a subfamily from Theorem 38 of \cite{IKMZ}.
By the uniqueness of tropical limit we have $C=\pi_1^{trop}(\widetilde{C})$.
\end{proof}
\begin{definition}
We refer to the map $\pi_2^{trop}|_{\widetilde{C}}:\widetilde{C}\rightarrow \TP^{1}$
as a tropical limit of the family $\gamma_t$, 
and we call it a \em tropical Gauss map.
\end{definition}
\begin{remark}
This construction extends in the case of hypersurfaces of higher dimensions.
However in this text, we will focus on the case of curves.
\end{remark}
\begin{remark}
A priori, the tropical curve $\widetilde{C}$ not only depends on $C$ but also on an approximation $(V_t)_{t\in A}$ of $C$, and even on the choice of the subfamily $A'\subset A$.
Nevertheless, the tropical limit of the ramification locus of $\gamma_t$
is a well-defined subset of $C$
and does not depend on these choices
as shown in this paper.
\end{remark}
We now define the degree of a projection $\pi:\Gamma \subset\R^{n}\mapsto \R^{k}$, where $\Gamma$ is a tropical curve, and where we identify $\R^{n}$ with $\R^{k}\times\R^{n-k}$. If the image $\pi(\Gamma)$ is a point, we say that $\pi$ is of degree $0$. Assume that $\pi(\Gamma)$ is not reduce to a point. Then if $x\in \pi(\Gamma)$ is a generic point, the fiber $\pi^{-1}(x)$ intersects $\Gamma$ transversally (meaning that any point of intersection is in the relative interior of an edge of $\Gamma$) at finitely many points. Denote by $e$ the edge of $\pi(\Gamma)$ containing $x$ and denote by $\Z_e$ the sublattice of $\Z^{k}\times\Z^{n-k}$ generated by $0\times \Z^{n-k}$ and $(\overrightarrow{v},0)$, where $\overrightarrow{v}$ is a primitive vector contained in $e$. Let $z\in \Gamma\cap \pi^{-1}(x)$, and denote by $\widetilde{e}$ the edge of $\Gamma$ containing $z$. Define the \textit{local tropical intersection number}
 $\iota_z (\Gamma,\pi^{-1}(x))$ as the index of the sublattice in $\Z_e$ generated by $\Z^{n-k}$ and a primitive vector contained in $\widetilde{e}$ (multiplied by the weight of $\widetilde{e}$). The degree of $\pi$ is defined as the sum 
 $$
 \deg(\pi):=\sum_{z\in \pi^{-1}(x)\cap\Gamma}\iota_z (\Gamma,\pi^{-1}(x)).
 $$
It follows from the balancing condition that the degree is independent of the point $x$ (see also Lemma 40 in \cite{BIMS}). 
The following proposition follows directly from the proof of Proposition $42$
of ~\cite{IKMZ} (the degree of $\pi_1$ is $1$ and the degree of $\pi_2$
is $2\mathrm{Area}(\Delta)$).
\begin{proposition}
One has $\deg(\pi_1^{trop})=1$ and $\deg(\pi_2^{trop})=2\mathrm{Area}(\Delta)$.
\end{proposition}

One can define similarly the local degree of $\pi$. Let $s\in\Gamma_1$ and let $U_s$ be a small neighborhood of $s$ only containing points from edges adjacent to $s$. If $\pi$ maps the neighborhood $U_s$ to a point, we say that $\pi$ is of local degree $0$ at $s$. If not, let $x\in \pi(\Gamma)\cap \pi(U_s)$ be a generic point. Define the local degree of $\pi$ at $s$ as the sum of local tropical intersection number $\iota_z(\Gamma,\pi^{-1}(x))$ over all $z\in U_s$. 


\section{First properties of tropical Gauss map}
\begin{proposition}
\label{technicallemma}
Let $\{V_t\subset\ctor\}_{t\in A}$ be a family of complex curves parameterised by the scaled sequence
$\alpha:A\to\R$. Assume that $\{V_t\subset\ctor\}_{t\in A}$ tropically converges to a smooth plane tropical curve $C$.
Let $e$ be an edge of $C$ of slope $p\in\QP^1$, and let $x$ belongs to the relative interior of $e$. If $x_t\in V_t$ satisfies $\lim_{t\in A}\mathrm{Log}_{\alpha (t)}(x_t)=x$, then
$$
\lim_{t\in A}\gamma_t (x_t)=p.
$$
\end{proposition}
\begin{proof}
Denote by $\{f_t=\sum_{j,k\in\Delta\cap\Z^2}a_{jk,t}z^jw^k\}_{t\in A}$ a family of polynomials defining the family of complex curves $\left(V_t\right)_{t\in A}$. Consider the tropical curve $C_x$ obtained by translating the tropical curve $C$ along the vector $-x\in\R^2$. Denote by $\{f_{t,x}\}_{t\in A}$ the family of polynomials given by
$$
f_{t,x}(z,w)=\alpha(t)^{-b}f_t(\alpha(t)^{x_1}z,\alpha(t)^{x_2}w)=\alpha(t)^{-b} \sum a_{jk,t}z^jw^k\alpha(t)^{b_{j,k}},
$$
where $x=(x_1,x_2)$ and $b_{j,k}=jx_1+kx_2$, and 
$$
b=b_{j_1,k_1}+a_{j_1,k_1}^{trop}=b_{j_2,k_2}+a_{j_2,k_2}^{trop},
$$ 
where $(j_1,k_1)$, and $(j_2,k_2)$ are the integer points on the edge (of the subdivision of $\Delta$) dual to the edge $e$.

The family of complex curves $\{V_{t,x}\}_{t\in A}=\{f_{t,x}=0\}_{t\in A}$ converges tropically to $C_x$. Furthermore, denoting by $\Phi_x$ the isomorphism 
$$
\begin{array}{ccc}
 V_{t,x} & \overset{\Phi_x }{\longrightarrow} & V_t \\
(z,w) & \mapsto & (\alpha(t)^{x_1}z,\alpha(t)^{x_2}w)
\end{array}
$$
one has the following commutative diagram

$$
\xymatrix{
    V_{t,x} \ar[r]^{\Phi_x} \ar[d]_{\gamma_{t,x}} & V_t \ar[d]^{\gamma_t} \\
    \CP^1 \ar[r]^{\mathrm{Id}} & \CP^1
  }.
$$
It is then enough to prove the proposition for the point $(0,0)\in C_x$ and the family $V_{t,x}$.
Since $b>b_{j,k}+a_{j,k}^{trop}$ for all $(j,k)\in\left(\Delta\cap\Z^2\right)\setminus\{(j_1,k_1),(j_2,k_2)\}$, there exist real numbers $c_{j,k}$ such that $c_{j_1,k_1}=c_{j_2,k_2}=0$ and $c_{j,k}>0$ otherwise such that  
$$
f_{t,x}(z,w)=\sum_{(j,k)\in\Delta\cap\Z^2} a_{jk,t}\alpha(t)^{-a_{j,k}^{trop}-c_{j,k}}z^jw^k.
$$
Since the tropical curve $C$ is smooth, there exists a $\Z$-affine automorphism $G$ of the plane such that $G(j_1,k_1)=(0,0)$ and $G(j_2,k_2)=(1,0)$. Write $G=M+(v_1,v_2)$, where $M\in GL_2(\Z)$ and $(v_1,v_2)\in\Z^2$. 
The map $M$ give rise to a multiplicative automorphism $\Phi_M$ of $\ctor$. 
Define
$$
g_t(z,w)=z^{v_1}w^{v_2}f_{t,x}(\Phi_M(z,w)),
$$
and 
$$
W_t=\{g_t=0\}.
$$
One has $\Phi_M(W_t)=V_{t,x}$ and the family of complex curves $(W_t)_{t\in A}$ tropically converges to the tropical curve $D$ which is the image of $C$ by $^{t}M^{-1}$. The edge of slope $p$ containing $(0,0)$ is then mapped to an edge of slope $\infty$. Furthermore one has the following commutative diagram:
$$
\xymatrix{
    W_{t} \ar[r]^{\Phi_M} \ar[d]_{\gamma_{g_t}} & V_{t,x} \ar[d]^{\gamma_{t,x}} \\
    \CP^1 \ar[r]^{M^{-1}} & \CP^1
  }.
$$
One deduce that it is enough to prove the statement for the point $(0,0)$ and the edge of slope $\infty$ in the tropical limit of the family of complex curves $(W_t)_{t\in A}$ given by the family of polynomials $g_t(z,w)$. One has
$$
g_t(z,w)=a_{j_1k_1,t}\alpha(t)^{-a_{j_1k_1}^{trop}}+a_{j_2k_2,t}\alpha(t)^{-a_{j_2k_2}^{trop}}z+R_t(z,w),
$$  
where $R_t(z,w)$ is  a Laurent polynomial such that all coefficients $r_{jk,t}$ satisfies $r_{jk,t}=O(\alpha(t)^{c})$, where $c< 0$ is some real number. 
Let's divide $g_t$ by $a_{j_1k_1,t}\alpha(t)^{-a_{j_1k_1}^{trop}}$ and make the change of variables $z\mapsto a_{j_2k_2,t}^{-1}\alpha(t)^{a_{j_2k_2}^{trop}} z$. It does not affect the limiting tropical curve since the tropical limit of $a_{j_2k_2,t}\alpha(t)^{-a_{j_2k_2}^{trop}}$ is zero. 
We end with a polynomial $h_t$ of the form
$$
h_t(z,w)=1+z+S_t(z,w),
$$
where as above $S_t(z,w)$ is  a Laurent polynomial such that all coefficients $s_{jk,t}$ satisfies $s_{jk,t}=O(\alpha(t)^{c})$. 
One has then
$$
\gamma_{h_t}(z_t,w_t)=\left[z_t(1+\dfrac{\partial S_t}{\partial z}):w_t\dfrac{\partial S_t}{\partial w}\right].
$$
Let $x_t=(z_t,w_t)\in \{h_t=0\}$ such that $\lim_{t\rightarrow +\infty}\mathrm{Log}_{t}(x_t)=(0,0)$. Since $z_t=-1-S_t(z_t,w_t)$, one obtains that $\lim_{t\rightarrow +\infty}z_t=-1$. One has also $\lim_{t\rightarrow +\infty}(z_t\dfrac{\partial S_t}{\partial z})=\lim_{t\rightarrow +\infty}(w_t\dfrac{\partial S_t}{\partial w})=0$, which proves the proposition.
\end{proof}
Taking the logarithm with base $\alpha(t)$, we obtain the following corollary, 
where $1_\T=0$ is the neutral element for the tropical multiplication. 
\begin{corollary}
\label{corollary: technicallemma}
Let $\{V_t\subset\ctor\}_{t\in A}$ be a family of complex curves parametrised by the scaled sequence
$\alpha:A\to\R$ which tropically converges to a smooth plane tropical curve $C$. Let $e$ be an edge of $C$ of slope $p\in\QP^1$, $x$ be a point of the relative interior of $e$, and $x_t\in V_t$ be a sequence such that $\lim_{t\in A}\mathrm{Log}_{\alpha(t)}(x_t)=x$. Then up to passing to a subsequence one has
\begin{itemize}
\item if $p\neq 0$ or $\infty$, then $\lim_{t\in A}\mathrm{Log}_{\alpha(t)}(\gamma_t(x_t))=\left[1_\T:1_\T\right]\in\TP^1$,
\item if $p=0$, then $\lim_{t\in A}\mathrm{Log}_{\alpha(t)}(\gamma_t(x_t))=\left[a:1_\T\right]\in\TP^1$, for some $a\in\left[-\infty,0\right]$,
\item if $p=\infty$, then $\lim_{t\in A}\mathrm{Log}_{\alpha(t)}(\gamma_t(x_t))=\left[1_\T:a\right]\in\TP^1$, for some $a\in\left[-\infty,0\right]$.
\end{itemize}
\end{corollary}

Since the tropical map $\pi_1^{trop}$ is of degree $1$, one has a continuous inclusion $i_C:C\hookrightarrow \widetilde{C}$ such that $\pi_1^{trop}\circ i_C=id$.
\begin{proposition}
We have $\pi_2^{trop}(s)=\left[1_{\T}:1_{\T}\right]\in\TP^1$ for any $s\in i_C(\mathrm{Vert}(C))$.
\end{proposition}
\begin{proof}
At least one of the edges adjacent to $s$ must have the slope different from $0$ and $\infty$. The proposition follows then immediately from Corollary \ref{corollary: technicallemma}.
\end{proof}

\begin{proposition}
\label{localdegree1}
For any $s\in i_C(\mathrm{Vert}(C))$, the tropical map $\pi_2^{trop}$ is of local degree $1$ at $s$.
\end{proposition}
\begin{proof}
Since the logarithmic Gauss map commutes with multiplicative translations, one can assume that $s=i_C(0,0)$. Performing a multiplicative automorphism $\Phi_M$ as in the proof of Proposition \ref{technicallemma}, denote by $(W_t)_{t\in A}$ the family of complex curves such that $\Phi_M (W_t)=V_t$. Denote by 
$\widetilde{W}_t$ the graph of the logarithmic Gauss map on $W_t$ and by $\Phi:=(\Phi_M,{}^t M^{-1})$ the map sending $\widetilde{W}_t$ to $\widetilde{V}_t$. Denote by $\widetilde{D}$ the tropical limit (up to passing to a subsequence) of $\widetilde{W}_t$, by $D=\pi_1^{trop}(\widetilde{D})$ and by $s'=i_D(0,0)$. By definition, the map $\Phi_M$ sends $\mathrm{Log}_{\alpha(t)}^{-1}(0,0)\subset W_t$ to $\mathrm{Log}_{\alpha(t)}^{-1}(0,0)\subset V_t$. Since $\widetilde{W}_t$ (resp., $\widetilde{V}_t$) is a graph over $W_t$ (resp., $V_t$), there exist a small neighborhood $U_{s'}$ of $s'$ and a small neighborhood $U_s$ of $s$ such that $\Phi(\mathrm{Log}_{\alpha(t)}^{-1}(U_{s'}))\subset\mathrm{Log}_{\alpha(t)}^{-1}(U_{s})$.
However, let $y\in \pi_2^{trop}(U_s)$ be a generic point, and let $(y_t)_{t\in A}\in\CP^1$ such that $\lim_{t\in A}\mathrm{Log}_{\alpha(t)}(y_t)=y$. It follows from Proposition 43 of \cite{BIMS} that the local degree of $\pi_2^{trop}$ at $s$ is equal to the intersection of $\widetilde{V}_t$ and $\pi_2^{-1}(y_t)$ in $\mathrm{Log}_{\alpha(t)}^{-1}(U_s)$, for sufficiently large $\alpha(t)$. Applying the map $\Phi$, one deduce that it is enough to prove the statement for the family of complex curves $(W_t)_{t\in A}$ given by the family of polynomials
$$
h_t(z,w)=1+z+w+S_t(z,w),
$$
where $S_t(z,w)$ is  a Laurent polynomial such that all coefficients $s_{jk,t}$ satisfies $s_{jk,t}=O(\alpha(t)^{c})$.
The graph of the logarithmic Gauss map over $W_t$ is then given by $$
\widetilde{W}_t=\left\lbrace z,w,\left[z(1+\dfrac{\partial S_t}{\partial z}):w(1+\dfrac{\partial S_t}{\partial w})\right]\mid (z,w)\in W_t\right\rbrace
$$
and the graph of the logarithmic Gauss map on the line $\mathcal{L}=\left\lbrace 1+z+w=0\right\rbrace$ is given by
$$
\widetilde{\mathcal{L}}=\left\lbrace z,w,\left[z:w\right]\mid (z,w)\in\mathcal{L}\right\rbrace.
$$ 
Denote by $\widetilde{L}$ the tropical limit of $\widetilde{\mathcal{L}}$. It follows from the description of $\widetilde{W}_t$ and $\widetilde{\mathcal{L}}$ that there exists a small neighborhood $U$ of the point $(0,0,0)$ in $\R^3$ such that $\widetilde{D}\cap U=\widetilde{L}\cap U$. But the map $\pi_2^{trop}:\widetilde{L}\rightarrow \TP^1$ is of degree one, which proves the proposition.
\end{proof}

\begin{example}
In the left hand side of Figure \ref{example: tropGaussmap}, we draw the tropical limit (when $t\rightarrow +\infty$) of $\widetilde{V}_t$, for $V_t=\{1+x+y+t^{-1}xy\}$. We draw also the image of this tropical limit under the first projection $\pi_1^{trop}$. In the right hand side of Figure \ref{example: tropGaussmap}, we did the same for $V_t=\{1+x+y+t^{-1}x^2\}$.
\begin{figure}[h]
\begin{tabular}{cc}
\includegraphics[height=8cm, width=4.5cm]{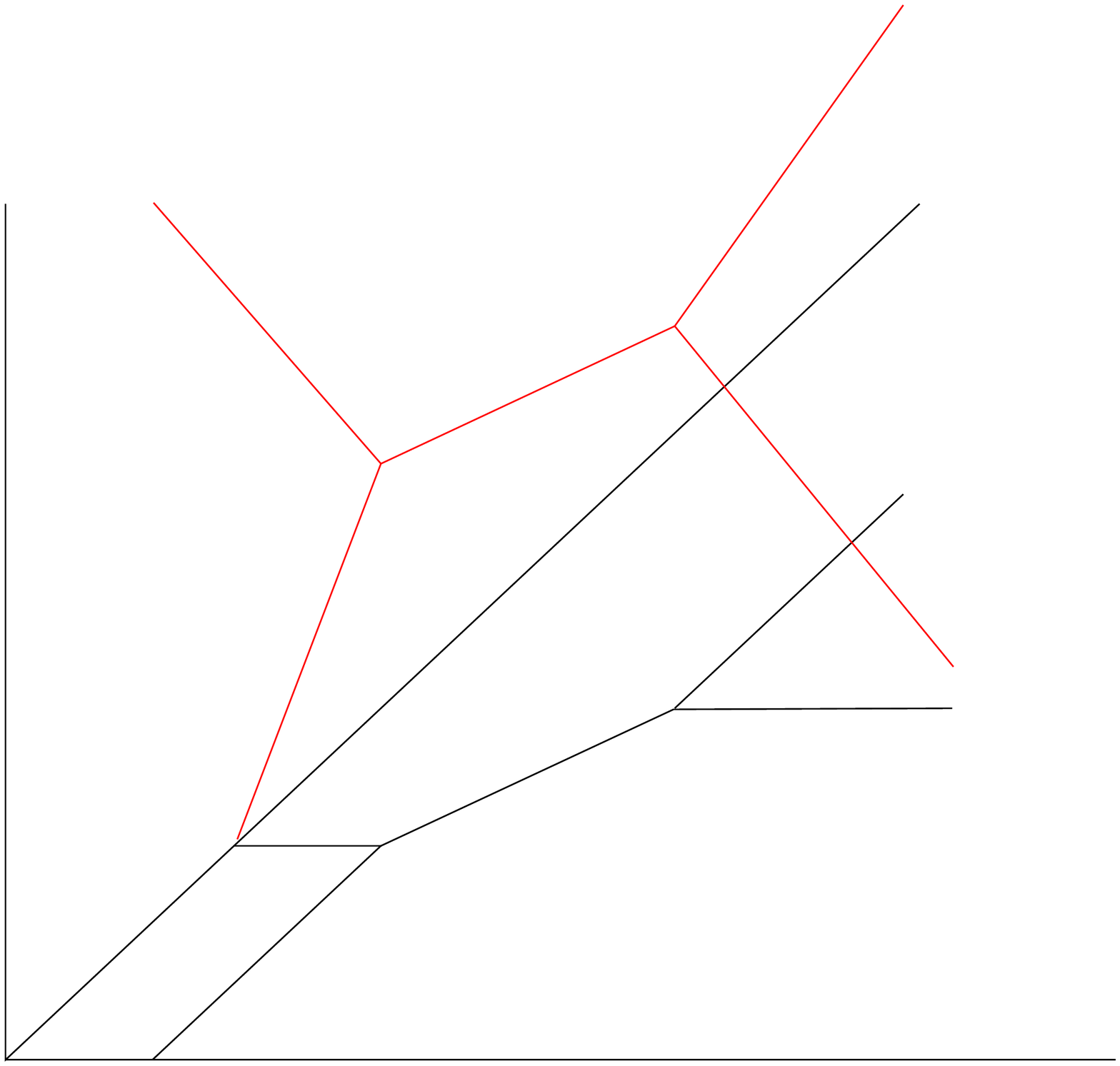}
  \put(-10,50){\tiny{$\T^2$}}
  \put(-150,170){\tiny{$\T^1\subset \TP^1$}}
  \put(-160,35){\tiny{$(-\infty,-\infty,-\infty)$}} & \hspace{7ex}
\includegraphics[height=8cm, width=4.5cm]{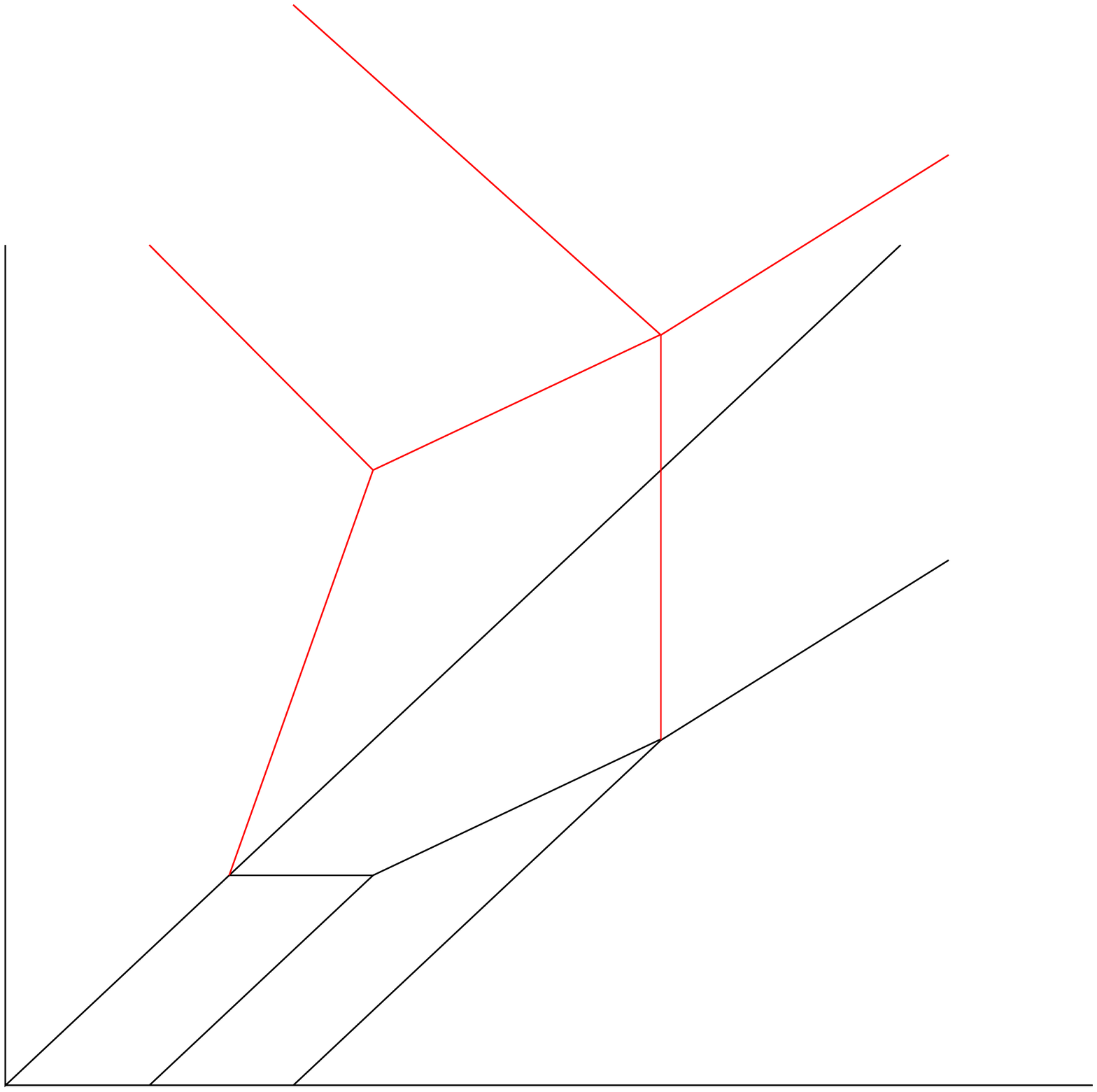}
\put(-10,50){\tiny{$\T^2$}}
  \put(-150,170){\tiny{$\T^1\subset\TP^1$}}
  \put(-160,35){\tiny{$(-\infty,-\infty,-\infty)$}}
\end{tabular}
\caption{Tropical limit of $\widetilde{V}_t$ for $V_t=\{1+x+y+t^{-1}xy\}$ and for $V_t=\{1+x+y+t^{-1}x^2\}$.}
\label{example: tropGaussmap}
\end{figure} 
\end{example}

\section{Proof of Theorem \ref{parabolicforcurves}, Proposition \ref{prop:realinflectionpts} and Corollary \ref{corollary:twistededges}} 
\label{proofmaintheorem}
The tropical Gauss map $\pi_2^{trop}$ contracts a priori some edges of $i_C(C)$. It is then not immediately clear how to locate the tropical limit of the logarithmic inflection points on $\widetilde{C}$. To overcome this difficulty, we embed $\widetilde{C}$ and $\TP^1$ in higher dimensional spaces such that the new tropical Gauss map does not contract any edges coming from $C$ (see Proposition \ref{prop: finiteatvertices}).  Then we prove that the tropical Gauss map is of local degree $2$ at each point coming from a midpoint of $C$ (see Proposition \ref{lemma: finite}).
\\Denote by $m_1,\cdots, m_l$ the midpoints of bounded edges of $C$, and by $x_t^1,\cdots, x_t^l$ some points in $V_t$ such that $\lim_{t\in A}\mathrm{Log}_{\alpha(t)}(x_t^i)=m_i$, for $1\leq i\leq l$. 
Put $\gamma_t(x_t^i)=\left[u_t^i:v_t^i\right]$, for $1\leq i \leq l$ and denote by 
$$
\begin{array}{cccc}
\Phi_t: & \CP^1 & \rightarrow & \CP^{l+1} \\
& \left[u:v\right] & \mapsto & \left[u:v:\varphi_t^1(u,v):\cdots :\varphi_t^l(u,v)\right]
\end{array},
$$
where $\varphi_t^i (u,v)=uv_t^i-vu_t^i$ are linear maps.

Consider the line $\mathcal{L}_t=\Phi_t(\CP^1)$, and denote by $L$ its tropical limit in $\TP^{l+1}$ (after passing to a subsequence), which exists by compactness theorem (see \cite{IKMZ}, Section $3.4$). Introduce also
$$
\widehat{V}_t=\left\lbrace (z_t,\Phi_t\circ\gamma_t(z_t)) \mid z_t\in V_t \right\rbrace\subset \ctor\times\CP^{l+1},
$$
and denote by $\widehat{C}\subset \R^2\times\TP^{l+1}$ its tropical limit (after passing to a subsequence). One has isomorphisms $\widehat{V}_t\simeq \widetilde{V}_t$ and $\mathcal{L}_t\simeq\CP^1$ given by projections. 
The projection on the second factor $\R^2\times \TP^{l+1}\rightarrow \TP^{l+1}$ restricts to a projection $\pi^{trop}:\hat{C}\rightarrow L$ of degree 
$2\mathrm{Area}(\Delta)$, the projection $\R^2\times \TP^{l+1}\rightarrow\R^2\times \TP^1$ restricts to a projection $p:\hat{C}\rightarrow \widetilde{C}$ of degree $1$ and the projection $\TP^{l+1}\rightarrow\TP^1$ restricts to a projection $q:L\rightarrow\TP^1$ of degree $1$. Furthermore, the following diagram commutes.

$$
\begin{CD}
\widehat{C}     @>\pi^{trop} >>  L\\
@VVpV        @VVqV\\
\widetilde{C}     @>\pi_2^{trop}>>  \TP^1\\
@VV\pi_1^{trop}V \\
C
\end{CD}
$$
Denote by $i_p:C\rightarrow \widehat{C}$ the continuous inclusion that is right-inverse to the projection $\pi_1^{trop}\circ p$, and denote by $s_L$ the image in $L$ of $\left[1_\T:1_\T\right]\in\TP^1$ by the continuous inclusion $\TP^1\hookrightarrow L$ that is right-inverse to $q$. It follows from Proposition \ref{localdegree1} that for any $s\in \mathrm{Vert}(C)$, one has $\pi^{trop}(i_p(s))=s_L$. In fact, if $\pi^{trop}(i_p(s))\neq s_L$, a small neighborhood of $i_C(s)$ would be contracted to $\left[1_\T:1_\T\right]$ by $\pi_2^{trop}$, contradicting Proposition \ref{localdegree1}. Moreover, since $\pi_2^{trop}$ is locally of degree $1$ at $i_C(s)$ for any $s\in\mathrm{Vert}(C)$ and since the projections $p$ and $q$ are both of degree $1$, we deduce that the tropical map $\pi^{trop}$ is locally of degree $1$ at $i_p(s)$ for any $s\in\mathrm{Vert}(C)$.

\begin{proposition}
\label{prop: finiteatvertices}
 Let $s$ be a vertex of $C$ and let $\widehat{e}$ be a bounded edge of $i_p(C)$ adjacent to $i_p(s)$. Then $\widehat{e}$ is not contracted by $\pi^{trop}$. 
\end{proposition}
\begin{proof}
Assume that $\widehat{e}$ is contracted by $\pi^{trop}$, and denote by $e$ the edge of $C$ containing $\pi_1^{trop}\circ p(\widehat{e})$. Let us first prove that $i_p(e)=\widehat{e}$. 
\\Assume that $\widehat{e}\varsubsetneq i_p(e)$. Then $\widehat{e}$ is adjacent to a vertex $s'\in Vert(\widehat{C})\setminus i_p(Vert(C))$. Since $\widehat{e}$ is contracted by $\pi^{trop}$, one has $\pi^{trop}(s')=\pi^{trop}(s)=s_L$. 
Since by definition every vertex is at least $3$-valent, there exists an edge adjacent to $s'$ not belonging to $i_p(Edge(C))$. Such an edge is contracted by $\pi_1^{trop}\circ p$ and so is not contracted by $\pi^{trop}$. Then the local degree of $\pi^{trop}$ at $s'$ is at least $1$. Since every vertex in $i_p(Vert(C))$ already contributes to $1$ to the degree of $\pi^{trop}$ and since $\#(Vert(C))=2\mathrm{Area}(\Delta)$, we conclude that the degree of $\pi^{trop}$ is at least $2\mathrm{Area}(\Delta)+1$, which is impossible. Then one has $i_p(e)=\widehat{e}$. 
\\It follows from the construction of $\widehat{C}$ that for any bounded edge $e_C$ of $C$, the set $i_p(e_C)$ consists of at least two edges (if $m_C$ is the middle point of $e_C$, there is a point of $\widehat{C}$ in a boundary divisor of $\TP^{l+1}$ projecting to $m_C$). Then $i_p(e_C)\neq \widehat{e}$ and $\widehat{e}$ cannot be contracted by $\pi^{trop}$.  
\end{proof}
We call a tropical map $f$ \textit{finite} over a subset $F$ if it does not contract any edge in $F$. 
\begin{proposition}
Let $e$ be a bounded edge of $C$, and let $m$ be the midpoint of $e$. The map $\pi^{trop}$ is finite over $i_p(e)$ and is of local degree $1$ at any point of $i_p(e)$ except at $i_p(m)$, where it is of local degree $2$.
\label{lemma: finite}
\end{proposition}
\begin{proof}
Denote by $T$ the image of $i_p(e)$ by the map $\pi^{trop}$. The set $T$ is a tree since $T\subset L$. Let us first show that $s_L$ is a vertex of valence $1$ of $T$. Assume that the valence of $s_L$ is at least $2$ and denote by $a_1,a_2$ two edges of $T$ adjacent to $s_L$. Since $T$ is a tree, the sets $a_1\setminus s_L$ and $a_2\setminus s_L$ are in two different connected components of $T\setminus s_L$. Since $i_p(e)$ is connected, it means that there are at least $3$ vertices of $i_p(e)$ mapped to $s_L$. As in the proof of Proposition \ref{prop: finiteatvertices} this give a contradiction with the degree of $\pi^{trop}$. Denote by $e_L$ the edge of $T$ adjacent to $s_L$. 
\\Assume first that the edges $e_1^1$ and $e_1^2$ of $i_p(e)$ adjacent to the vertices over $s_L$ join over the edge $e_L$, see Figure \ref{figure1}.
\begin{figure}[h]
\centering
\includegraphics[height=6cm,width=12cm]{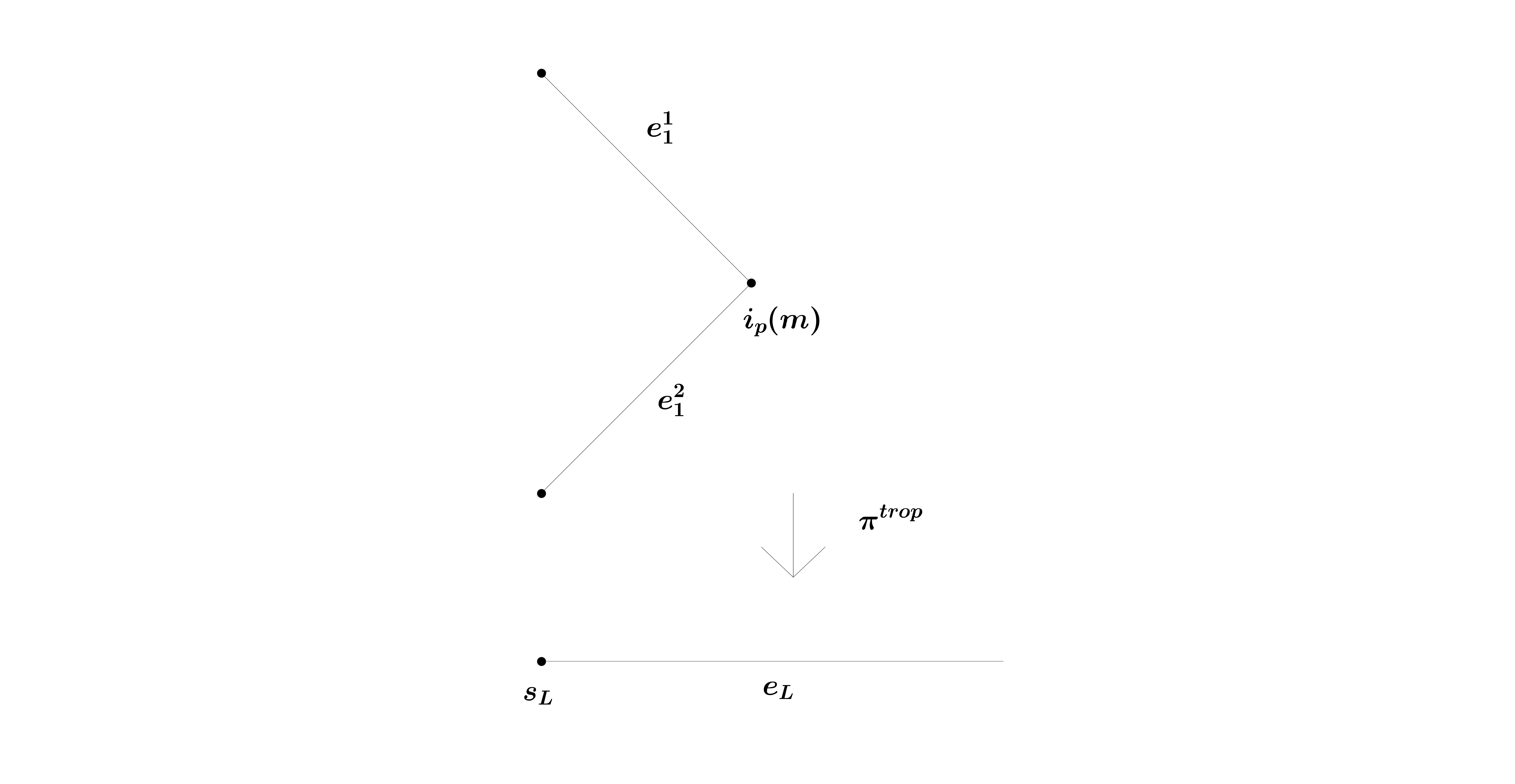}
\caption{ \label{figure1} The two edges $e_1^1$ and $e_1^2$ join over the edge $e_L$.}
\end{figure}
\\In this case, one has $i_p(e)=e_1^1\cup e_1^2$. Since the local degree of $\pi^{trop}$ at the two vertices in $i_p(Vert(C))$ adjacent to $i_p(e)$ is equal to one, the local degree of $\pi^{trop}$  on the edges $e_1^1$ and $e_1^2$ is equal to one. Then the vertex adjacent to $e_1^1$ and $e_1^2$ is the vertex $i_p(m)$ (the lengths of $e_1^1$ and $e_1^2$ are equal). Moreover, the map $\pi^{trop}$ is of local degree $2$ at $i_p(m)$. In fact, if the local degree would be more than $2$, then by the balancing condition $s_L$ would have at least $3$ preimages by $\pi^{trop}$ over the set $(\pi_1^{trop}\circ p)^{-1}(e)$, which is impossible.
\\Assume now that the two edges $e_1^1$ and $e_1^2$ do not join over $e_L$, and denote by $s_1$ the other vertex of $T$ adjacent to $e_L$, see Figure \ref{figure2}. 
\begin{figure}[h]
\centering
\includegraphics[height=6cm,width=12cm]{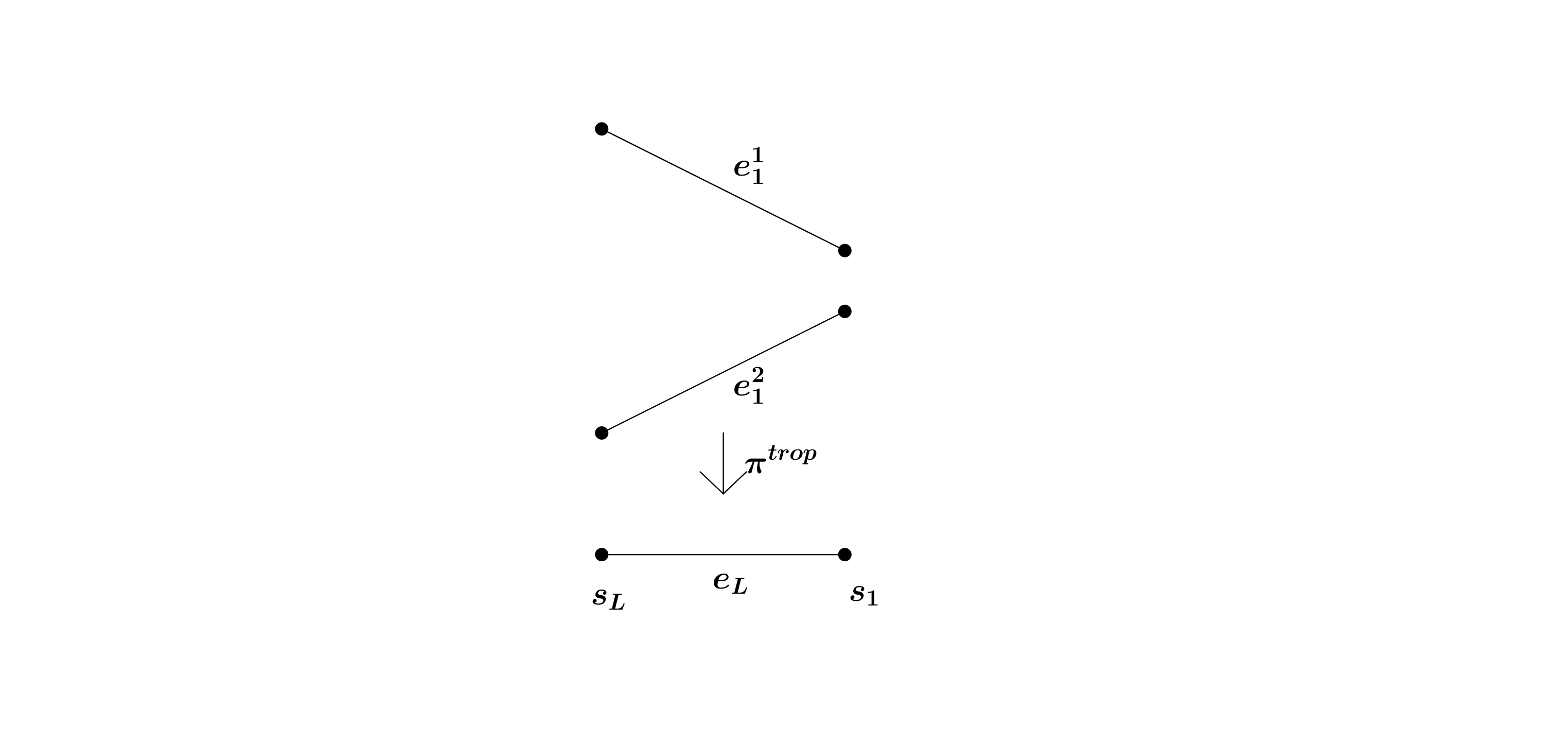}
\caption{\label{figure2}The two edges $e_1^1$ and $e_1^2$ do not join over the edge $e_L$.}
\end{figure}
\\ In this case, the preimage of $s_1$ consists of $2$ vertices at which the map $\pi^{trop}$ is of local degree one. As in the proof of Proposition \ref{prop: finiteatvertices}, one see that every edge of $i_e(C)$ adjacent to $i_e(s_1)$ is not contracted by $\pi^{trop}$. One conclude then easily by induction.

\end{proof} 
\begin{proof}[Proof of Theorem \ref{parabolicforcurves}]
Let $e$ be a bounded edge of $C$, and denote by $v_1$ and $v_2$ the vertices adjacent to $e$. Consider the fibration $\lambda_t:V_t\rightarrow C$ as defined in ~\cite{Mik04}. It follows from Proposition \ref{localdegree1} that there exists a small neighborhood $\mathcal{N}_i$ of $v_i$, $i=1,2$ in $C$ such that for $t$ big enough, the map $\gamma_t\vert_{\lambda_t^{-1}(\mathcal{N}_i)}$ is of degree one on its image, see Figure \ref{retraction2}.
\begin{figure}[h]
\centering
\includegraphics[height=4cm,width=6cm]{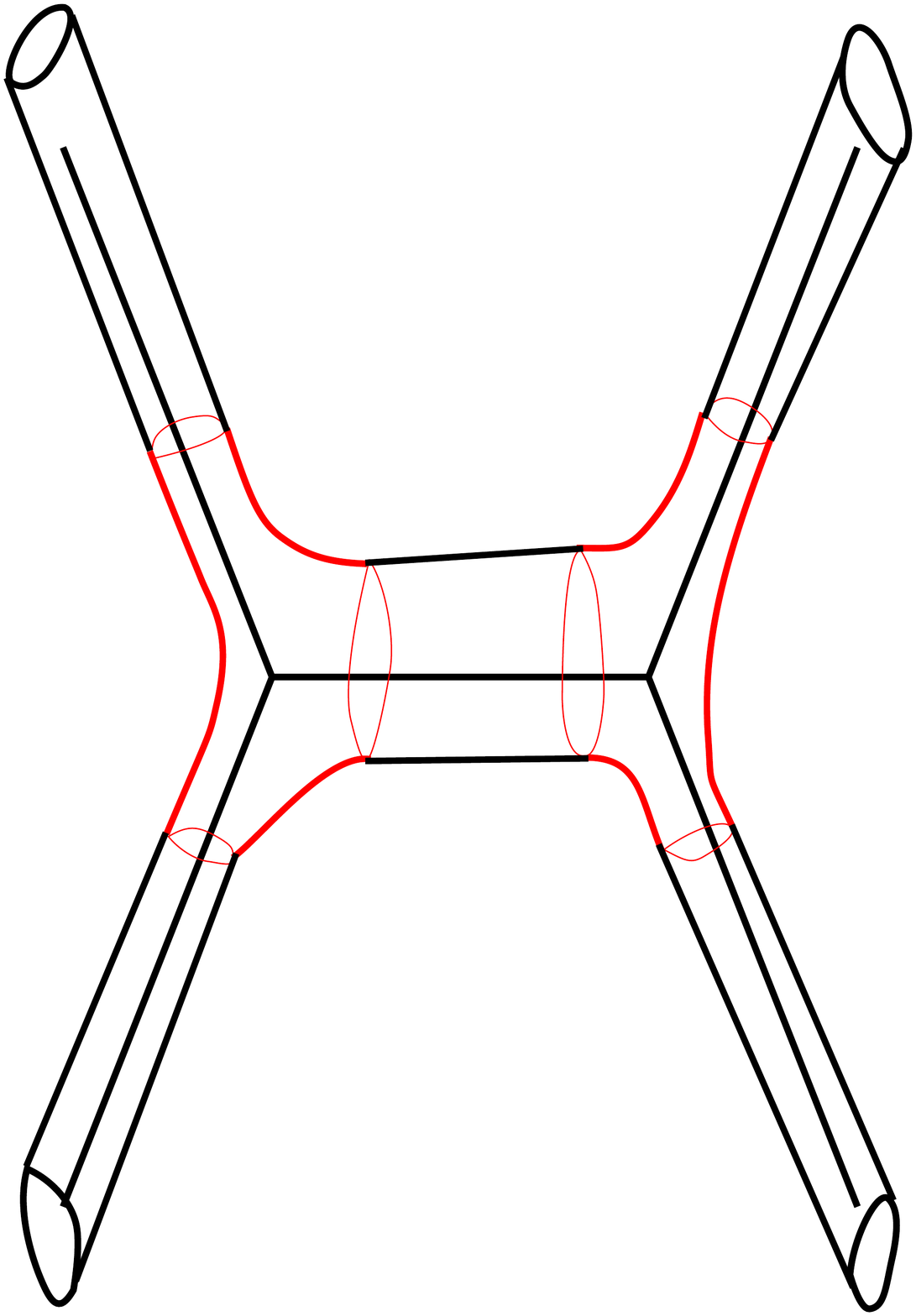}
\put(-85, 57){\tiny{$e$}}
\put(-120, 50){\tiny{$v_1$}}
\put(-60,50){\tiny{$v_2$}}
\caption{\label{retraction2} The sets $\lambda_t^{-1}(\mathcal{N}_i)$ in red inside $V_t$.}
\end{figure}
 It follows from Proposition \ref{technicallemma} that the image by $\gamma_t$ of a boundary component of $\lambda_t^{-1}(\mathcal{N}_i)$ converges to the slope of the edge supporting this boundary component. Consider a small circle $\Gamma\subset\CP^1$ centered at the slope of $e$ and of radius $\varepsilon$. The image by $\gamma_t$ of the boundary component of $\lambda_t^{-1}(\mathcal{N}_i)$ supported above $e$ is contained in the same connected component of $\CP^1\setminus\Gamma$ as the slope $s_e$ of $e$. Then for $t$ big enough, $\Gamma$ is contained in $\gamma_t\left(\lambda_t^{-1}(\mathcal{N}_i)\right)$, and the set $\gamma_t^{-1}\vert_{\lambda_t^{-1}(\mathcal{N}_i)}(\Gamma)$ is a circle $\Gamma_i\subset \lambda_t^{-1}(\mathcal{N}_i)$ (see Figure \ref{retraction3}).
\begin{figure}[h]
\centering
\includegraphics[height=4cm,width=6cm]{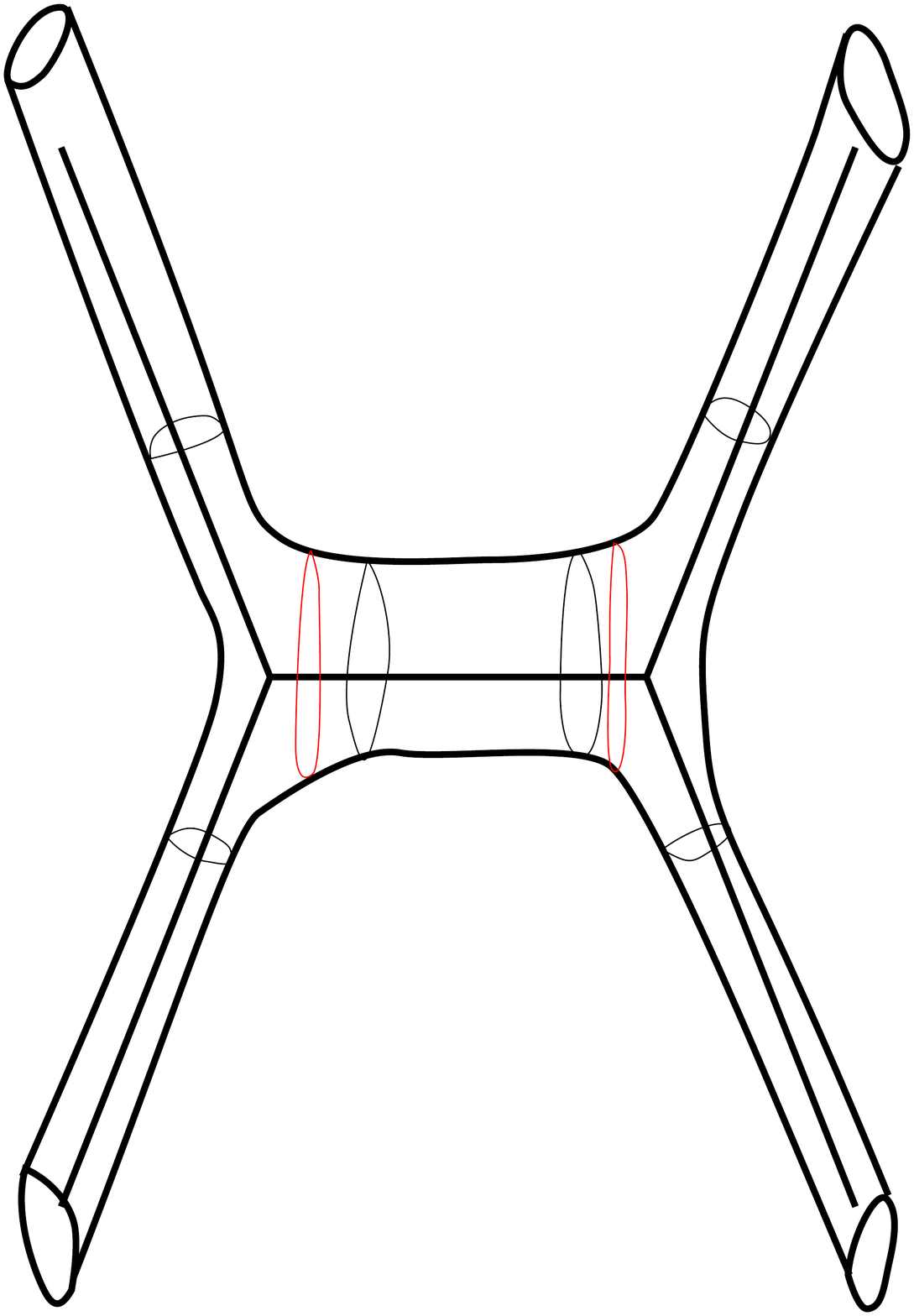}
\put(-119, 36){\tiny{$\Gamma_1$}}
\put(-61,36){\tiny{$\Gamma_2$}}
\caption{\label{retraction3} The circles $\Gamma_i$ in red inside the sets $\lambda_t^{-1}(\mathcal{N}_i)$.}
\end{figure}
 Note that one of the two connected components of $\lambda_t^{-1}(\mathcal{N}_i)\setminus\Gamma_i$ only contains the boundary component of $\lambda_t^{-1}(\mathcal{N}_i)$ supported on $e$. If not, then there would exist a path in $V_t\setminus\Gamma_i$ with ends $a_t$ and $b_t$, such that the tropical limit of $a_t$ lie in the interior of $e$ and the tropical limit of $b_t$ lie in the interior of another edge adjacent to $v_i$ (see Figure \ref{retraction4}).
 \begin{figure}[h]
\centering
\includegraphics[height=4cm,width=6cm]{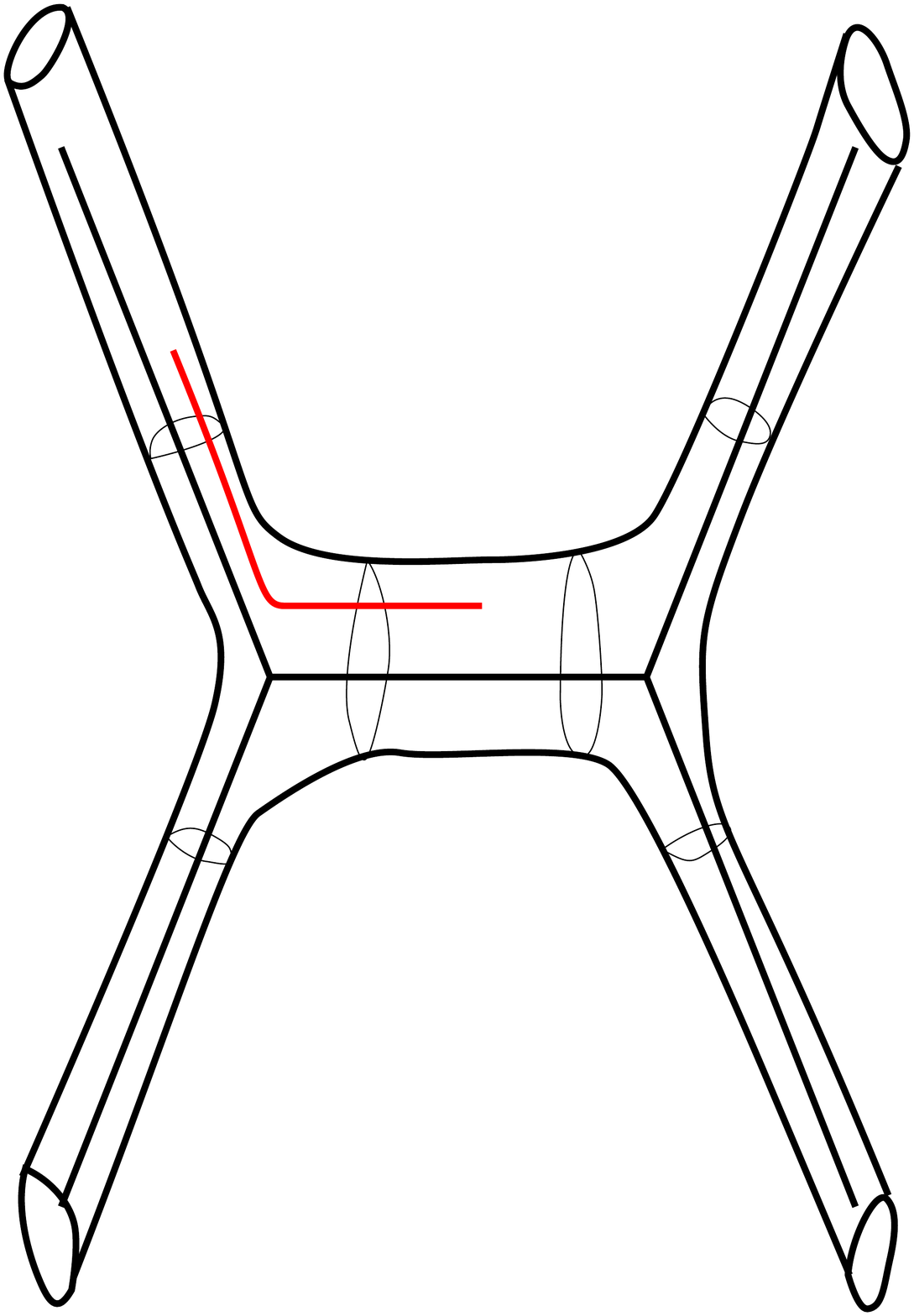}
\put(-145,85){\tiny{$b_t$}}
\put(-80,57){\tiny{$a_t$}}
\caption{\label{retraction4} The path with ends $a_t$ and $b_t$ inside $V_t$.}
\end{figure}
But it follows from Proposition \ref{technicallemma} that for $t$ big enough, such a path should intersect $\Gamma_i$.   
\\ The restriction of $\gamma_t$ to the annulus $\mathcal{A}$ supported on $e$ with boundary $\Gamma_1\cup\Gamma_2$ gives us a map of degree $2$ to the disc containing $s_e$ with boundary $\Gamma$. In fact, the preimage of $\Gamma$ by $\gamma_t$ in $\mathcal{A}$ consists exactly of $\Gamma_1\cup\Gamma_2$ since the map $\gamma_t$ is of degree $2\mathrm{Area}(\Delta)$ and since from Proposition \ref{localdegree1} for any vertex $v$ of $C$ there exists a small neighborhood $\mathcal{N}$ of $v$ such that the map $\gamma_t$ restricted to $\lambda_t^{-1}(\mathcal{N})$ is of degree $1$. By the Riemann-Hurwitz formula, such a map has two critical points. The tropical limit of a scaled sequence of critical points must be a point with the local degree of $\pi^{trop}$ greater than one. Proposition \ref{lemma: finite} now implies the theorem as $\rho V_t$ are the critical points
of $\gamma_t$.
\end{proof}

\begin{proof}[Proof of Proposition \ref{prop:realinflectionpts}]
Suppose that $\sigma_{v_1^E}\sigma_{v_2^E}=-1$.
Then there is a branch of $\R V_t$
corresponding to $E$ in the positive quadrant. The tropical
limit of the branch of $\R V_t$ corresponding to $E$ in the
positive quadrant is non-convex (see Figure \ref{picture:inflection}) if and only
if $\sigma_{v_3^E}\sigma_{v_4^E}=-1$. 
\begin{figure}[h]
\includegraphics[height=2cm, width=2.5cm]{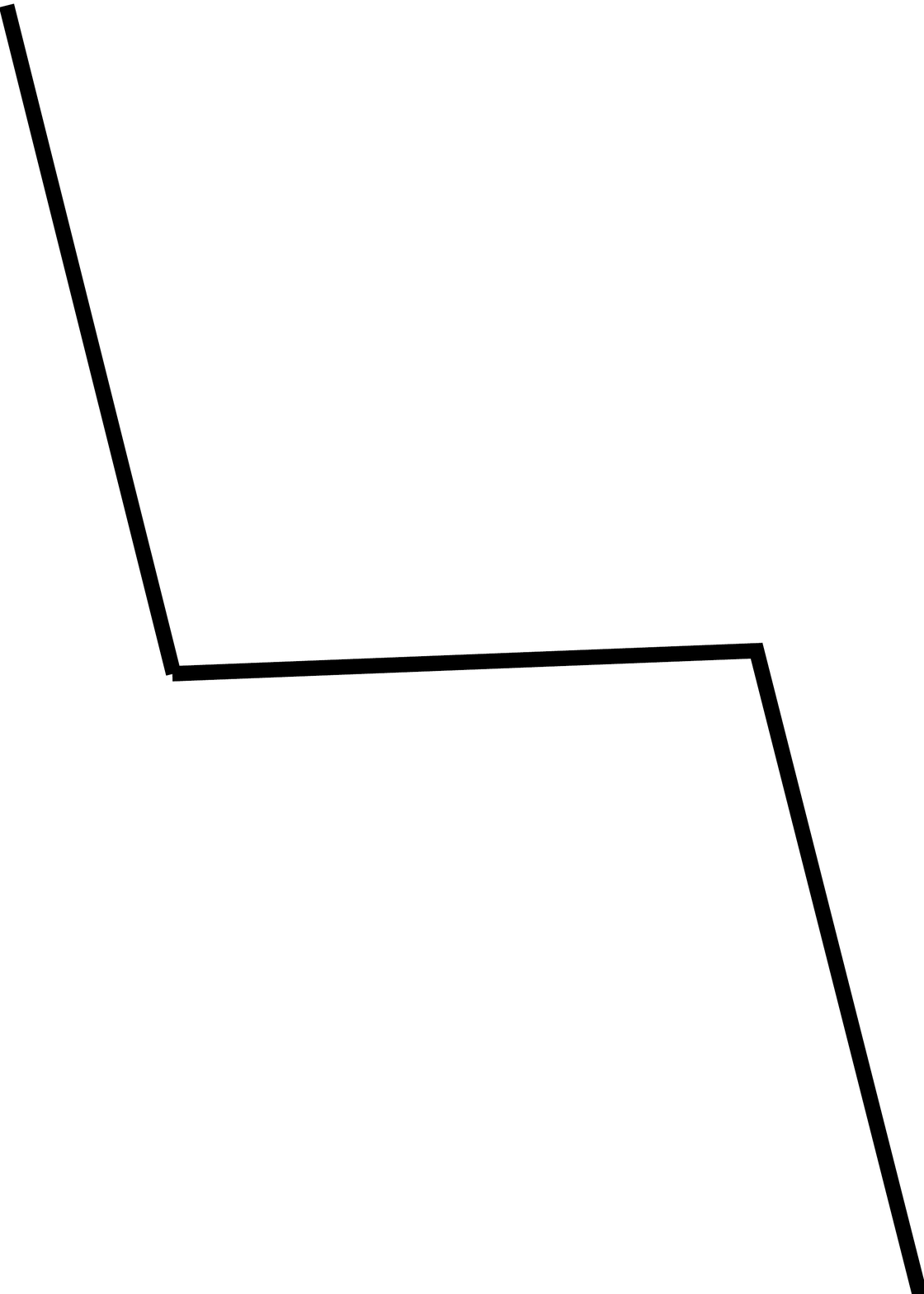}
\put(-35,30){\tiny{$E$}}
\caption{\label{picture:inflection} Non-convex tropical limit of a branch of $\R V_t$ corresponding to an edge $E$.}
\end{figure} 
Thus both log-inflection
points corresponding to $E$ must be real (and belong to
two different quadrants).

If $\sigma_{v_3^E}\sigma_{v_4^E}=1$ then the tropical limit of both corresponding branches are convex (see Figure \ref{picture:convex}). 
\begin{figure}[h]
\includegraphics[height=2cm, width=2.5cm]{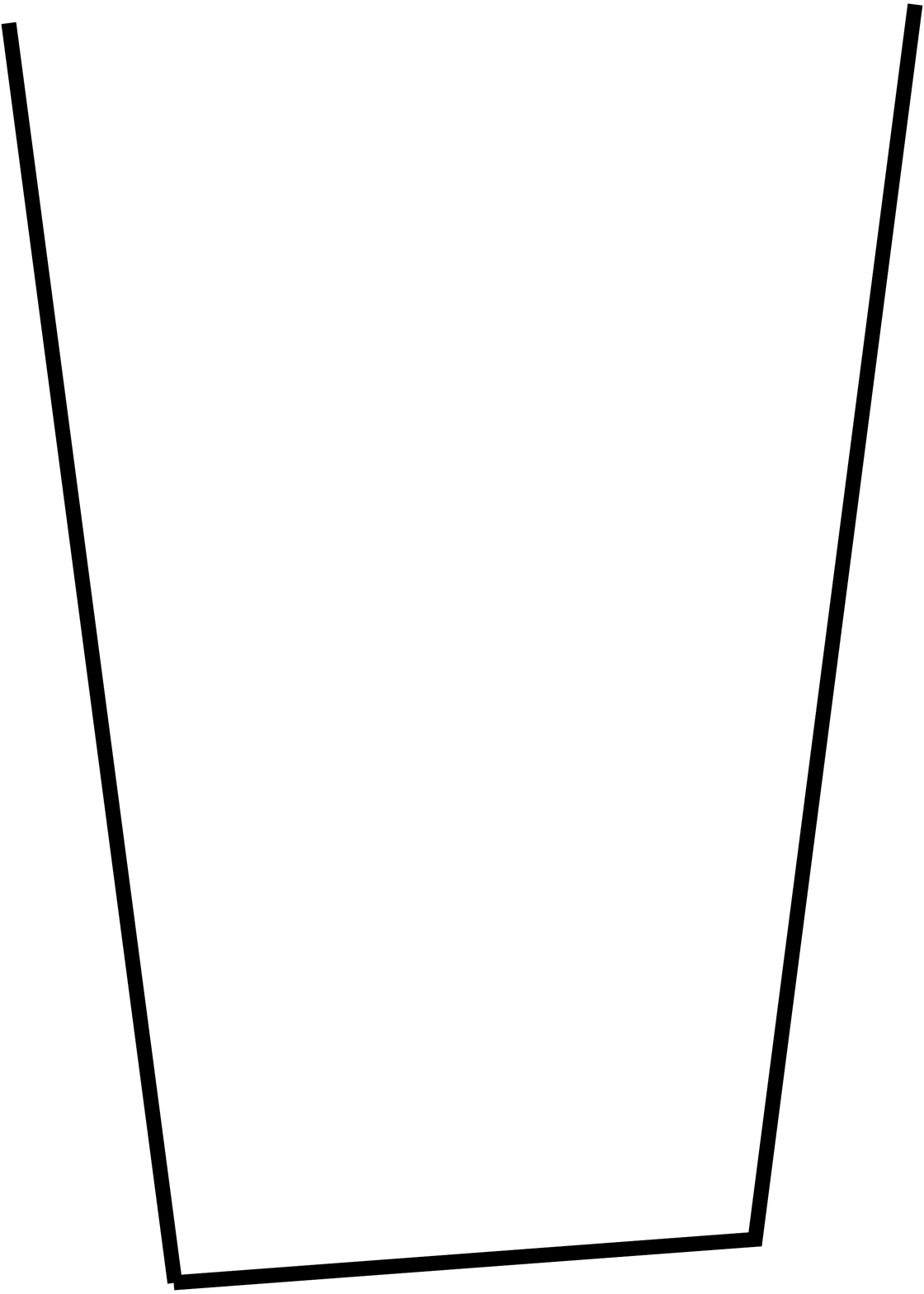}
\put(-35,5){\tiny{$E$}}
\caption{\label{picture:convex} Convex tropical limit of a branch of $\R V_t$ corresponding to an edge $E$.}
\end{figure}
Suppose that one of the branches $\mathcal{B}$ of $\R V_t$ corresponding to $E$ contains two real log-inflection points. Consider
the second branch $\mathcal{B}'$
of $\R V_t$
corresponding to $E$.
If $\mathcal{B}'$ is also
logarithmically non-convex then
we have more than two inflection points for $E$,
which contradicts Theorem \ref{parabolicforcurves}.
If $\mathcal{B}'$ of $\R V_t$ is
logarithmically convex 
then we can find a tangent line to $\mathcal{B}'$ parallel to an inflection point of $\mathcal{B}$ (in the logarithmic
coordinates). This gives us at least three points (counted with multiplicity) with the same image
under the logarithmic Gauss map $\gamma_t$.
There are $2\mathrm{Area}(\Delta)-2$
vertices of $C$ not adjacent to $E$.
It follows from Proposition \ref{localdegree1}
that each of them
contributes one point of $V_t$ with the same image of $\gamma_t$.
We arrive to a contradiction with $\deg \gamma_t=2\mathrm{Area}(\Delta)$.

This finishes the proof of the proposition in the case $\sigma_{v_1^E}\sigma_{v_2^E}=-1$.
The case $\sigma_{v_1^E}\sigma_{v_2^E}=1$ is obtained by applying
the coordinate change $z\mapsto -z$ or/and $w\mapsto -w$.
\end{proof}

\bibliographystyle{alpha}
\bibliography{biblio}
\end{document}